\documentclass[a4paper,11pt,reqno]{amsart}

\usepackage{type1cm}        
\usepackage{graphicx}        
\usepackage{multicol}        
\usepackage[bottom]{footmisc}

\usepackage[utf8]{inputenc}
\usepackage{amsmath}
\usepackage{amssymb}
\usepackage{xcolor}
\usepackage{hyperref}
\usepackage{cleveref}

\usepackage{enumerate}
\usepackage{verbatim}
\usepackage{float}

\newcommand{\Ric}{\mathrm{Ric}}

\def\A{{\mathcal A}}

\def\k{{\mathfrak k}}

\newcommand{\BE}{{\sf BE}}
\newcommand{\RCD}{{\sf RCD}}
\newcommand{\CD}{{\sf CD}}
\newcommand{\MCP}{{\sf MCP}}

\def\N{{\mathbb N}}
\def\Z{{\mathbb Z}}
\def\R{{\mathbb R}}

\DeclareMathOperator{\vol}{\mathsf{vol}}

\newcommand{\m}{\mathfrak{m}}

\renewcommand{\d}{\,\mathrm{d}}
\newcommand{\ap}[1][t,s]{\hat P_{#1}}
\renewcommand{\liminf}{\varliminf}
\renewcommand{\limsup}{\varlimsup}
\newcommand{\Hess}{\mathrm{Hess}}

\newcommand{\rcd}{\mathrm{RCD}}
\newcommand{\ent}{\mathrm{Ent}}

\newcommand{\p}{\mathcal{P}}
\newcommand{\spt}{\mathrm{spt}\,}

\numberwithin{equation}{section}

\theoremstyle{plain}
\newtheorem{thm}{Theorem}[section]
\newtheorem{theorem}[thm]{Theorem}
\newtheorem{proposition}[thm]{Proposition}
\newtheorem{definition}[thm]{Definition}

\newtheorem{corollary}[thm]{Corollary}

\newtheorem{remark}[thm]{Remark}
\newtheorem{example}[thm]{Example}

\begin{document}

\title{Synthetic approaches to Ricci flows}

\author{Matthias Erbar$^\dagger$}
\author{Marco Flaim$^*$}
\author{Eric Hupp$^*$}
\author{Zhenhao Li$^\dagger$}
\author{Timo Schultz$^\dagger$$^\ddag$$^*$}
\author{Karl-Theodor Sturm$^*$}

\address{$^\dagger$Faculty of Mathematics \\university
         Bielefeld University \\
         Postfach 10 01 31 \\
         33501 Bielefeld \\
         Germany.}
\address{$^\ddag$Department of Mathematics and Statistics\\
        University of Jyv\"askyl\"a \\
        P.O.Box 35 (MaD)\\
        FI-40014\\
        Finland.}
\address{$^*$Institut f\"ur angewandte Mathematik\\
  Universit\"at Bonn\\
   Endenicher Allee 60\\
   53115 Bonn\\
   Germany.}

\email{matthias.erbar@math.uni-bielefeld.de}  
\email{flaim@iam.uni-bonn.de}
\email{hupp@iam.uni-bonn.de}      
\email{zhenhao.li@math.uni-bielefeld.de}
\email{timo.m.schultz@jyu.fi}
\email{sturm@uni-bonn.de}

%
%

\begin{abstract}
We review different notions of synthetic Ricci flow that apply to time-dependent families of metric measure spaces and which are based on properties of the heat flow, ideas from optimal transport, and the asymptotic behaviour of volumes. Each notion equivalently characterises (weighted) Ricci flow for smooth families of weighted Riemannian manifolds. We discuss the features of the different notions on various examples.
\end{abstract}

\maketitle

\section{Introduction}
\label{sec:intro}

A smooth manifold $M$ equipped with a smooth one-parameter family of
Riemannian metrics $(g_t)_{t\in I}$ evolves according to \emph{Ricci
  flow}, if
\begin{equation}
  \label{eq:RF}
  \partial_tg_t = -2\Ric_{g_t}\;.
\end{equation}
Since the groundbreaking work of Hamilton \cite{H82,H95}, Ricci flow has
received a lot of attention and has become a powerful tool in many applications, most prominently in Perelman's work on the Poincar\'e conjecture and Thurston's geometrisation conjecture \cite{P02,P03,P03b}, see also \cite{CZ06,KL08,MT07}. For a detailed account on the Ricci flow we refer e.g.~to \cite{CK04}.
\medskip

A challenging feature of Ricci flow is the development of singularities in finite time. Therefore it seems desirable to obtain robust characterisations of Ricci flow that make it possible to consider evolutions of non-smooth spaces. Among many recent contributions in this direction let us highlight the following:
Ricci flows with irregular or incomplete initial metrics \cite{Sim02, Sim12, KL12, GT13, MRS15, Yin10, ST21}; a weak notion of Ricci flow in 3 dimensions and construction of a canonical Ricci flow through singularities as the unique limit of Ricci flows with surgery \cite{KL17, KL20, BK22};  a compactness theory for super Ricci flows and partial regularity results for non-collapsed limits of Ricci flows \cite{Bam20,Bam23}; characterisations of Ricci flow in terms of functional inequalities on the path space equipped with the Wiener measure \cite{Haslhofer-Naber18, CT17}.
\medskip

In this paper we present several synthetic notions of Ricci flow that are based on optimal transport, properties of the heat flow, ideas from optimal transport, and the asymptotic behaviour of volumes. This builds on a sequence of developments connecting Ricci curvature and optimal transport for both static and time-dependent manifolds and metric measure spaces, that we briefly recall.

Lower bounds for the Ricci curvature of Riemannian manifolds $(M,g)$ can be characterised in terms of semiconvexity properties of entropy-like functionals on the Wasserstein space $\mathcal P_2(M)$ \cite{Cordero-Erausquin2001,vonRenesseSturm}. This gives rise to a powerful and stable definition of synthetic lower Ricci bounds for metric measure spaces (mm-spaces for short) \cite{sturm2006--1,Sturm06-2,LottVillani09}. This Lagrangian formulation of  synthetic lower Ricci bounds indeed turns out to be equivalent to the Eulerian formulation as proposed by Bakry and \'Emery much earlier in the setting of Markovian semigroups \cite{BakryEmeryDiffHypercontr, AGS3, EKS}.
 
Synthetic characterisations of upper Ricci bounds are much more subtle. Several, not necessarily equivalent such characterisations are proposed in \cite{sturm2017remarks}. Instead of uniform estimates they rely on asymptotic estimates. In general, they are less stable.
   
The dynamic counterpart to manifolds of nonnegative Ricci curvature are super-Ricci flows of Riemannian manifolds. They also allow for a broad variety of synthetic characterisations, both in Eulerian and Lagrangian descriptions, and a far reaching extension to super-Ricci flows of mm-spaces \cite{McCann-Topping,Sturm2018Super,Kopfer-Sturm2018,KopferSturm21}. The concept of super-Ricci flows, however, has limited applications. Deeper insights and stronger results for flows of manifolds and mm-spaces are available only  for Ricci flows or, more generally,  for flows satisfying a parabolic curvature condition as introduced in \cite{muller2010monotone} and analysed recently in \cite{flaim}.\medskip 
  
In Section \ref{sec:weakrough}, we provide synthetic definitions Ricci flow in terms of short time-asymptotics of contractivity estimates for the heat flow and convexity/concavity properties of the entropy along optimal transports. To this end we complement the notions of synthetic notion of super-Ricci flow mentioned above with suitable notions of sub-Ricci flow which are the dynamic counterparts to the synthetic notions of upper Ricci bounds for mm-spaces. The corresponding results have been obtained in \cite{ELS25}.

In Section \ref{sec:minimal}, we present a different perspective on the challenge of synthetic characterisations of Ricci flows that is pursued in the forthcoming paper \cite{FHS}. Here, we seek to identity \emph{minimal} super-Ricci flows in terms of the asymptotic behaviour of the volume of small balls. In the static case, the central idea is that in oder to characterise Ricci flatness of a closed manifold, it suffices to characterise non-negativity of the Ricci curvature and non-positivity of the integrated scalar curvature. The latter is captured in the  asymptotic behaviour of volume.

\section{Weak and rough Ricci flows}
\label{sec:weakrough}

In this part we present two synthetic notions of Ricci flow that are based on Wasserstein contraction properties of the heat flow on one hand and on convexity properties of the entropy along optimal transports on the other hand.

\subsection{Background and Heuristics} 

We start by recalling the charcaterisation of lower Ricci bounds using optimal transport and the heat flow obtained \cite{Cordero-Erausquin2001, vonRenesseSturm}.

Let $(M,g)$ be a Riemannian manifold. We denote by $W$ the $L^2$ Wasserstein distance on the space $\mathcal P_2(M)$ of probability measures with finite second moment over $M$ built from the Riemannian distance $d$ associated to $g$. For a probability measure $\mu$ on $M$ the Boltzmann entropy is given by 
\[\ent(\mu)=\int \rho\log \rho\; \d \vol_g\;,\]
provided $\mu=\rho\vol_g$ is absolutely continuous w.r.t. the volume measure $\vol_g$, and by $\ent(\mu)=+\infty$ else. We denote by $P_t=e^{t\Delta}$ the heat semigroup generated by the Laplace-Beltrami operator $\Delta$ and by $\hat P_t$ 
the dual semigroup acting on measures. Now, the following are equivalent

\begin{enumerate}
\item[($0$)] $\Ric_g\geq 0$\;;
\item[($1$)] \emph{geodesic convexity of the entropy}: for any constant speed Wasserstein geodesic $(\mu^a)_{a\in[0,1]}$ and all $a\in[0,1]$ we have
\[\ent(\mu^a)\leq (1-a)\ent(\mu^0) + a \ent(\mu^1)\;; \]
\item[($2$)] \emph{Wasserstein contractivity of the heat flow}: for all $\mu,\nu\in \mathcal P_2(M)$
\[W(\hat P_t\mu, \hat P_t\nu) \leq W(\mu,\nu) \;.\]
\end{enumerate}
The latter properties are robust and can be used to give a synthetic definition of Ricci curvature bounds for non-smooth spaces with the only structure required being a distance and a reference measure. Starting from the pioneering works by Sturm \cite{sturm2006--1,Sturm06-2} and Lott and Villani \cite{LottVillani09} a rich and still rapidly growing theory of metric measure spaces with synthetic Ricci bounds has been developed. For a concise partial overview, we refer e.g. to \cite{Amb18}.

McCann and Topping \cite{McCann-Topping} and later Sturm \cite{Sturm2018Super} showed that this approach can be generalised to a dynamic setting of a time-dependent family of metrics $(g_t)$ on the manifold $M$ to obtain the following characterisation:

\begin{enumerate}
\item[($0_{\sf dyn}$)] $(M,g_t)$ is a super-Ricci flow, i.e. $\partial_t g_t \geq -2\Ric_{g_t}$\;;
\item[($1_{\sf dyn}$)] \emph{dynamic convexity of the entropy}: for all $t$ and any constant speed geodesic $(\mu^a)_{a\in[0,1]}$ in $(P(M),W_t)$ we have
\[\partial_a\big|_{a=1-}\ent_t(\mu^{a})-\partial_a\big|_{a=0+}\ent_t(\mu^{a})\geq -\frac12 \partial_t W_t(\mu^0,\mu^1)^2\;; \]
\item[($2_{\sf dyn}$)] \emph{Wasserstein contractivity of the heat flow}: for any $s\leq t$ and $\mu,\nu\in \mathcal P(M)$
\[W_s(\hat P_{t,s}\mu, \hat P_{t,s}\nu) \leq W_t(\mu,\nu) \;.\]
\end{enumerate}
Here, $ P_{t,s}$ denotes the heat propagator and $\hat P_{t,s}$ its dual, i.e. $P_{t,s}f$ gives the solution at time $t$ of $\partial_t u = \Delta_t u$ with initial datum $u=f$ at time $s$. Further $\Delta_t$, $W_t$, and $\ent_t$ denote the Laplace-Beltrami operator, the Wasserstein distance, and the entropy associated with the metric tensor $g_t$. 
Sturm \cite{Sturm2018Super} and Kopfer and Sturm \cite{Kopfer-Sturm2018} have used this to define a synthetic notion of super-Ricci flow for time-dependent families of metric measure spaces $\left(X, d_t, m_t\right)_{t\in I}$. Suitable regularity assumptions are needed for the second notion in order to ensure existence of the heat flow.
\medskip

We complement the above notions of super-Ricci flow with corresponding notions of sub-Ricci flow. This will be achieved by reversing the inequalities in ($1_{\sf dyn}$) and ($2_{\sf dyn}$) up to an arbitrarily small error for sufficiently localised transports and small times. We build on recent ideas and results in the static case \cite{sturm2017remarks} to encode Ricci upper bounds. This has also been used in the Lorentzian setting \cite{Mondino-suhr} to give characterisations of the Einstein equations in terms of optimal transport. 

We briefly give some intuition behind the idea of reversing the inequalities in ($1_{\sf dyn}$) and ($2_{\sf dyn}$).  A well-known formal computation in optimal transport \cite{vonRenesseSturm} shows that the second derivative of the entropy along a Wasserstein geodesic $(\mu^a)$ w.r.t.~$\m_t$ is given by
\[
\frac{d^2}{da^2}\ent_t(\mu^a) = \int_M\Big[\Ric_{g_t}(\nabla\psi^a) +\|\Hess\;\psi^a\|^2_{\rm HS}\Big]{\rm d}\mu^a\;,
\]
where $(\psi^a)$ is the family of Kantorovich potentials associated with the geodesic $(\mu^a)$. Assuming super-Ricci flow $\partial_tg_t \geq -2\Ric_{g_t}$ and neglecting the positive term $\|\Hess \psi^a\|_{\rm HS}^2$ formally yields the dynamic convexity inequality in $(1_{\sf dyn})$. However, to obtain a concavity estimate under the assumption of sub-Ricci flow $\partial_tg_t \leq -2\Ric_{g_t}$, the Hessian term can not be neglected as it has the wrong sign. 
The key idea is to consider transports that are sufficiently concentrated around a single geodesic on $M$ and where $\psi^a$ thus behaves almost linearly, so that the Hessian term becomes an arbitrarily small error. This leads to the almost concavity estimate in \eqref{eq:EntSub-intro} below. A similar reasoning applies to the non-expansion of the heat flow under super-Ricci flow and the almost non-contraction under sub-Ricci flow for sufficiently ,,linear" transports.

\subsection{Synthetic concepts and consistency}\smallskip

In the following, let $(X,d_t,\m_t)_{t\in I}$ be a time-dependent family of metric measure spaces with $I\subset \R$ an interval. That is, $X$ is a Polish space, and $(d_t)$, $(\m_t)$ are Borel families of distance functions inducing the given topology on $X$ and locally finite Radon measures on $X$. Let $W_t$ denote the $L^2$-Wasserstein distance on the space $\p_t(X)$ of probability measures with finite second moment w.r.t. $d_t$ and $\ent_t$ denote the relative entropy w.r.t. $\m_t$, i.e. for a probability measure $\mu\in \p_t(X)$ we set
\[\ent_t(\mu)=\int\rho\log\rho \d \m_t\;,\]
provided $\mu=\rho \m_t$ is absolutely continuous and $\ent_t(\mu)=+\infty$ else. 
\medskip

The first notion of Ricci flow we consider is based on dynamic convexity and almost concavity of the entropy.\medskip

\begin{definition}
The family $(X,d_t,\m_t)_{t\in I}$  is called a {\bf weak Ricci flow} if it satisfies
\begin{enumerate}[(i)]
\item (\emph{weak super-Ricci flow}): the entropy is strongly dynamically convex, i.e. for a.e. $t\in I$ and every $W_t$-geodesic $(\mu^a)_{a\in [0,1]}$ in $\p_t(X)$ with finite entropy at endpoints, the function $a\mapsto \ent_t(\mu^{a})$ is absolutely continuous on $[0,1]$ and
\begin{align}\label{ineq:weaksuperRF-intro}
    \partial^+_a \ent_t(\mu^a)\lvert_{a=1-}-\partial^-_a \ent_t(\mu^a)\lvert_{a=0+}\ge -\frac12\partial^-_tW_{t^-}^2(\mu^0,\mu^1),
\end{align}
holds, and

\item (\emph{weak sub-Ricci flow}) for a.e. $t$ and every $\varepsilon>0$, there exists an open cover $\{U_i\}$ such that for all $i$ and every open subsets $V_0,V_1\subset U_i$, there exists a $W_t$--Wasserstein geodesic $(\mu^a)$ with $\spt\mu^0\subset V_0$, $\spt\mu^1\subset V_1$, and
\begin{align}\label{eq:EntSub-intro}
    \partial_a^+\ent_t(\mu^a)\lvert_{a=1-}-\partial_a^-\ent_t(\mu^a)\lvert_{a=0+}\le -\frac12\partial^-_t W_{t^-}^2(\mu^0,\mu^1)+\varepsilon W_t^2(\mu^0,\mu^1).
\end{align}
\end{enumerate}
\end{definition}

Here, $\partial^\pm_a\big|_{a=c\mp}$ denotes the upper/lower right/left derivative.
\medskip

The second notion of Ricci flow we consider is based on expansion properties of the heat flow. Let us assume that $(X,d_t,\m_t)_{t\in I}$ satisfies additional regularity properties, namely, we require a uniform log-Lipschitz control in the time parameter $t$ on the distance functions $d_t$ and the measures $\m_t$ and that for each fixed $t$, the space $(X,d_t,\m_t)$ satisfies the Riemannian curvature-dimension condition RCD$(K,N)$ for some $K\in \R$ and $N\in[1,\infty)$. Under these conditions, Kopfer and Sturm \cite{Kopfer-Sturm2018} have shown the existence of a (dual) heat flow $(P_{t,s})$.\medskip

\begin{definition}
We call $(X,d_t,\m_t)_{t\in I}$  a {\bf rough Ricci flow} if it satisfies
\begin{enumerate}[(i)]
    \item (\emph{rough super-Ricci flow}) for all $s\leq t$ and all $x,y\in X$
    \begin{align}\label{ineq:defroughsuperRicci-intro}
        W^2_s(\ap{\delta_x},\ap{\delta_y})\le d^2_t(x,y)\;;
    \end{align}
    \item (\emph{rough sub-Ricci flow}) for a.e. $t$ and every $\varepsilon>0$ there exists an open cover $\{U_i\}_{i\in\N}$ such that for every $i$ and all $x,y\in U_i $, there exists $s_0<t$ so that
    \begin{align}\label{ineq:defroughsubRicci-intro}
        W^2_s(\ap{\delta_x},\ap{\delta_y})\ge d^2_t(x,y)-\varepsilon d^2_t(x,y)(t-s)
    \end{align}
    for every $s\in (s_0,t)$.
\end{enumerate}
\end{definition}

Note that the notions of weak/rough sub-Ricci flow above both formalise the idea of locally reverting the inequalities in ($1_{\sf dyn}$) and ($2_{\sf dyn}$) respectively up to a small error. It will further be convenient to consider the following localised quantities. 

For $\varepsilon>0$ and $t\in I$, $x,y\in X$ we set
 \begin{align*}
  &\eta_{\varepsilon}(t,x,y):= \\
  &\quad\inf\Big\{
\frac1{W_t^2(\mu^0,\mu^1)}\cdot\Big[ \partial^{+}_a \ent_t(\mu^a)\big|_{a=1}-\partial_a^{-} \ent_t(\mu^a)\big|_{a=0}
+\frac12\partial_{t}^{-}W_{t^-}^2(\mu^0,\mu^1)\Big]\Big\}\;,
\end{align*}
where the infimum is taken over all $W_t$-geodesics 
$(\mu^a)_{a\in[0,1]}$ such that $\mu^0$ and $\mu^1$ have finite entropy and are supported in balls of radius $\varepsilon$ w.r.t. $d_t$ around $x$ and $y$ respectively. Moreover, we set
\[
\eta(t,x,y):= \sup_{\varepsilon> 0}\eta_{\varepsilon}(t,x,y)\;,\quad
  \eta^*(t,x):=\limsup_{y,z\to x}\eta(t,y,z)\;.
\]
Bounds on these quantities describe the dynamic convexity/concavity of the entropy for transports between measures concentrated around the points $x,y$. Similarly, we define for $t\in I$ and $x,y \in X$
\begin{equation}\label{eq:theta-intro}
    \vartheta(t,x,y):= -\liminf_{s\nearrow t}\frac1{t-s}\log\frac{W_s(\ap\delta_x,\ap\delta_y)}{d_t(x,y)}\;,
\end{equation}
as well as
\begin{equation}\label{eq:thetastar-intro}
    \vartheta^*(t,x)=\limsup_{y,z\to x}\vartheta(t,y,z)\;,
\end{equation}
describing the short time asymptotics of the transport cost between two heat flows starting from $x,y$ respectively.\medskip

For families of smooth weighted Riemannian manifolds the notions of weak/rough Ricci flow indeed yield a characterisation of the classical notion of (weighted) Ricci flow. Let $(M,g_t)_{t\in I}$ be a smooth family of closed Riemannian manifolds with Riemannian distance $d_t$ and let $(f_t)_{t\in I}$ be a smooth family of functions on $M$. The associated metric measure spaces $(M,d_t, e^{-f_t}\vol_{g_t})_{t\in I}$ will be called a \emph{smooth flow}. The weighted Ricci tensor is defined by
\[\Ric_{f_t}:= \Ric_{g_t}+\Hess f_t\;.\]

\begin{theorem}\label{thm:consistency-intro}
Let $(M,d_t, e^{-f_t}\vol_{g_t})_{t\in I}$ be a smooth flow. Then the following are equivalent:
\begin{enumerate}[(i)]
\item The family of metric measure spaces $(M, d_t,e^{-f_t}\vol_{g_t})$ is a weak Ricci flow;
\item The family of metric measure spaces $(M, d_t,e^{-f_t}\vol_{g_t})$ is a rough Ricci flow;
\item for all $t\in I$ and $x,y\in M$ we have
\[\eta(t,x,y) \geq 0\;,\quad \eta^*(t,x)\leq 0\;;\]
\item for all $t\in I$ and $x,y\in M$ we have
\[\vartheta(t,x,y) \geq 0\;,\quad \vartheta^*(t,x)\leq 0\;;\]
\item $(M,g_t,e^{-f_t}\vol_{g_t})$ is weighted Ricci flow, i.e.
\[\partial_t g_t = - 2\Ric_{f_t}\;.\]
\end{enumerate}
\end{theorem}

\medskip

The characterisation above is obtained through a detailed local analysis of the short-time asymptotics of transport costs along the heat flow and of the dynamic convexity/almost concavity of the entropy along geodesics. Let us define the Ricci flow excess of a smooth flow given for $t\in I$ and $x,y\in M$ by
\[
{\rm RFex}(t,x,y) :=\inf \frac{1}{d_t(x,y)^2}\int_0^1\big[\Ric_{f_t}(\dot\gamma^a) +\frac12 \partial_tg_t(\dot\gamma^a)\big] {\rm d}a\;,
\]
where the infimum is taken over all geodesics from $x$ to $y$. We have the following estimates:
\begin{theorem} We have for all $t\in I$, $x,y\in M$:
\begin{align}\label{eq:excess-intro1}
{\rm RFex}(t,x,y) \leq \vartheta(t,x,y)\;.
\end{align}
For every $t_0\in I$ there is $\varepsilon>0$ such that for all $t\in I$, $x,y\in M$ non-conjugate with $|t-t_0|, d_t(x,y)<\varepsilon$: 
\begin{align}\label{eq:excess-intro2}
 \vartheta(t,x,y) \leq {\rm RFex}(t,x,y) + \sigma_t \tan^2\big(\sqrt{\sigma_t} d_t(x,y)\big)\;,
\end{align}
where $\sigma_t$ is an upper bound on the modulus of the Riemann tensor along the geodesic from $x$ to $y$. The same estimates hold for $\eta$ in place of $\vartheta$.
\end{theorem}
\medskip

Note that the weight $e^{-f_t}$ on the volume measure presents an additional degree of freedom that influences the evolution of the metric in a weighted Ricci flow $\partial_tg_t=\Ric_{f_t}$ but whose own evolution is not constraint. It is therefore desirable to be able to single out \emph{unweighted} Ricci flows through a synthetic characterisation. This can be achieved by a dimensional refinement of the notion of weak/rough super-Ricci flow proposed in \cite{Sturm2018Super,Kopfer-Sturm2018}. 

\begin{definition}
A family of m.m.s. $\left(X,d_t,\m_t\right)_{t\in I}$ is called a
\begin{enumerate}[(i)]
\item  \emph{weak $N$-super-Ricci flow} for $N\in [1,\infty]$ if the inequality \eqref{eq:EntSub-intro} above is strengthened to 
\begin{align*}
 &\partial^+_a \ent_t(\mu^a)\lvert_{a=1-}-\partial^-_a \ent_t(\mu^a)\lvert_{a=0+}\\
 \ge &-\frac12\partial^-_tW_{t^-}^2(\mu^0,\mu^1) +\frac{1}{N} \Big|\ent_t(\mu^1)- \ent_t(\mu^0)\Big|^2\;.
\end{align*}
\item \emph{rough $N$-super-Ricci flow} if for all $s<t$ and $\mu,\nu\in \mathcal P_t(X)$ we have
\[
        W_s^2(\ap\mu,\ap\nu)\le W_t^2(\mu,\nu)-\frac2N\int_{[s,t]}[\ent_r(\hat P_{t,r}\mu)-\ent_r(\hat P_{t,r}\nu)]^2\d r\;.
\]
\end{enumerate}
\end{definition}

Now, a smooth flow $(M,g_t,e^{-f_t}\vol_{g_t})$ of $n$-dimensional manifolds is a weak/rough $N$-super-Ricci flow, if and only if we have
\[\partial_tg_t\geq -2 \Ric_{N,f_t}\quad \text{with}\quad\Ric_{N,f_t}:= \Ric_{g_t} +\Hess f_t -\frac{1}{N-n}\nabla f_t\otimes \nabla f_t\;,\]
where the latter is the so-called weighted-$N$-Ricci tensor. We then obtain the following synthetic characterisation of Ricci flows. 
\begin{corollary}\label{cor:non-collapsed Ricciflow}
A smooth flow $(M,d_t, e^{-f_t}\vol_{g_t})_{t\in I}$ is a \emph{Ricci flow}, i.e. $\partial_tg_t=-2\Ric_{g_t}$ and $f_t$ is constant for all $t$, if and only if it is a weak/rough sub-Ricci flow and a weak/rough $N$-super-Ricci flow for some $N\in[1,\infty)$.
\end{corollary}
This is due to the observation that the combination of the bounds
\[ -\Ric_{f_t}\geq  \frac12\partial_tg_t\geq -\Ric_{N,f_t}= -\Ric_{f_t}+\frac{1}{N-n}\nabla f_t\otimes\nabla f_t\]
for some $N\geq n$ necessarily implies that $\nabla f_t\equiv 0$. Based on this result we call a family of m.m.s $\left(X,d_t,\m_t\right)_{t\in I}$ a \emph{non-collapsed weak/rough Ricci flow} if it is a weak/rough sub-Ricci flow and weak/rough $N$-super-Ricci flow for some finite $N$. We conjecture that such flows are indeed non-collapsed in the sense that the reference measure $\mathfrak m_t$ is a multiple of the Hausdorff measure w.r.t. $d_t$ for a.e. $t$.
\medskip

\subsection{Relating weak and rough Ricci flows}\smallskip

Let $(X,d_t,\m_t)_{t\in I}$ be a time-dependent family of m.m.s. with suitable regularity properties such that the heat flow is well-defined. From \cite{Kopfer-Sturm2018} it is known that the notions weak and rough super-Ricci flow are equivalent. 
On the other hand, weak and rough sub-Ricci flow turn out not to be equivalent. 
Indeed, for static spaces it has been shown in \cite{sturm2017remarks} that upper Ricci bounds in terms of transport cost asymptotics for the heat flow imply upper Ricci bounds in terms of almost concavity of the entropy. However, the latter notion does not detect the positive Ricci curvature in the vertex $o$ of a cone while the former does. More precisely, a Euclidean cone is Ricci flat in terms of convexity/almost concavity of the entropy, while $\vartheta^+(x,o)=+\infty$ for any point $x$ as shown in \cite{Erbar-Sturm}.
\medskip

We introduce a relaxation of the quantity $\vartheta^*$ 
given by
\begin{align*}
  \vartheta^\flat(t,x,y):= \lim_{\varepsilon\to 0}\inf\left\{-\limsup_{s\nearrow t}\frac1{t-s}\log\frac{W_s(\ap\mu^0,\ap\mu^1)}{W_t(\mu^0,\mu^1))}\right\}\;,
\end{align*}
where we infimise over all $\mu^0,\mu^1$ with $\spt(\mu^0)\subset B_{t}(x,\varepsilon)$ and $\spt(\mu^1)\subset B_t(y,\varepsilon)$, the balls of radius $\varepsilon$ around $x,y$ w.r.t. $d_t$. Trivially, we have $\vartheta^\flat\leq \vartheta^-$. In fact, $\vartheta^\flat$ is the lower semi-continuous envelope of $\vartheta^-$. Then we have the following result.

\begin{theorem}\label{thm:relation-intro}
For a.e. $t\in I$ we have $\vartheta^\flat(t,x)= \eta^*(t,x)$ for all $x\in X$. 
In particular, any rough sub-Ricci flow is also a weak sub-Ricci flow.
\end{theorem}

We will see below that the reverse implication fails, i.e. rough sub-Ricci flow is strictly stronger then weak sub-Ricci flow.

\subsection{Examples for weak/rough Ricci flows}

\subsubsection*{Smooth flows.} As already discussed, any smooth (sub/super)-Ricci flow $(M, d_t,e^{-f_t}\vol_{g_t})_{t\in I}$ is also a weak/rough (sub/super)-Ricci flow, this holds in particular for smooth Ricci flows starting from non-smooth initial data as considered e.g. in \cite{Sim12}, i.e $I=[t_0,t_1)$ and $(X,d_t,\m_t)$ approximates to a non-smooth m.m.s. as $t\searrow t_0$.

\subsubsection*{Static cones}

The constant flow of an $\rcd(0,N)$ metric measure space can be interpreted as a weak or rough $N$-super-Ricci flow. 
Notably, this includes any Euclidean $(N-1)$-cone over an $\rcd(N-2,N-1)$ metric measure space, provided $N \geq 2$.
However, by Theorem 1.1 in \cite{Erbar-Sturm}, the only Euclidean $N$-cone with rough Ricci curvature bounded above by 0, and satisfying $\rcd(K,N')$ for some $K, N' \in \mathbb{R}$, is Euclidean space $\mathbb{R}^{N+1}$. 
Therefore, the only static Euclidean $N$-cone that forms a rough Ricci flow is Euclidean space.

\subsubsection*{Spherical suspension}
Let $(X,d,\m)$ be a metric measure space.
The \emph{$N$-spherical suspension} is the mms defined on the set $\Sigma(X)=X\times [0,1]/\sim$ where $(x,r)\sim (y,r)$ $\Leftrightarrow$ $r=s=0$ or $r=s=\pi$ i.e. $\mathcal{S}=X\times\{0\}$ and $\mathcal{N}=X\times \{\pi\}$ are contracted to a point, the south and north pole respectively.
$\Sigma(X)$ is equipped with the following distance $d_{\Sigma}$ and measure $\m_{\Sigma}$: for $(x,s),(x',s')\in \Sigma(X)$
\begin{align}
\label{eq:d_conic}
    &\cos\left(d_{\Sigma}((x,s),(x',s')) \right):= \cos s\cos s'+\sin s\sin s'\cos(d(x,x')\wedge \pi
    )\;,\\
 &\d \m_{\Sigma}(x,s):=\d \m(x)\otimes \sin^N(s)\d s\;.
\end{align}
It is well known after Ketterer \cite{Ketterer-JMPA15} that when the base space $(X,d,\m)$ satisfies $\rcd(N-1,N)$ for some $N\geq 1$ and $\mathrm{diam}(X)\leq \pi$, then $(\Sigma(X),d_{\Sigma},\m_{\Sigma})$ satisfies $\rcd(N,N+1)$.
Consider the following time scaling of the metric
\begin{equation}
    d_t:=(1-2Nt)^{\frac12}d_{\Sigma},\quad \m_t:= c_t\m_{\Sigma},
\end{equation}
with $c:I\to \R$ a suitably regular function.
Then the time-dependent space $(\Sigma(X),d_t,\m_t)_{t\in(0,\frac{1}{2N})}$ is a $(N+1)$-super-Ricci flow, see e.g. \cite[Proposition 2.7]{Sturm2018Super} for a short explanation.

The following result shows that when the base space has a prescribed constant Ricci curvature, the time-scaled suspension becomes a weak Ricci flow.
However it qualifies as a rough sub-Ricci flow only in the absence of singularities at the poles, i.e., when it is a sphere.
\begin{theorem}\label{prop:suspension}
    Let $M$ be an $n$-dimensional Einstein manifold with $\Ric_{g_M}\equiv (n-1)g_M$.
    Then the time-scaled spherical suspension $(\Sigma(M),d_t,\m_t)_{t\in (0,\frac{1}{2n})}$ is a weak Ricci flow, where 
    \begin{equation}\label{eq:scaling}
        d_t:=(1-2nt)^{\frac12}d_{\Sigma},\quad \m_t:=(1-2nt)^{\frac{n+1}{2}}\m_{\Sigma}.
    \end{equation}
     Moreover, it is a rough super Ricci flow and it is rough sub-Ricci flow if and only if $M$ is the unit sphere $\mathbb{S}^{n}$ with the round metric and a multiple of volume measure.
\end{theorem}

That families of suspensions are weak Ricci flows essentially relies on the observation made in \cite{BacherSturm14} that optimal transports do not concentrate mass in the poles of the suspensions.

In particular, the scaled suspension $(\Sigma(M),d_t,\m_t)_{t\in [0,\frac{1}{8})}$ for $M=S^2(1/\sqrt{3})\times S^2(1/\sqrt{3})$ is not a rough Ricci flow, which was conjectured in \cite{Kopfer-Sturm2018}. The part on sub-Ricci flow in the result above also robust enough to apply to general $\rcd$-spaces. Let $\Sigma(X)$ be a spherical suspension over an $\rcd(n-1,n)$ space $(X,d,\m)$ for some $n\geq 1$. Then the scaled suspension $(\Sigma(X),d_t,\m_t)_{t\in [0,\frac{1}{2n})}$, given by \eqref{eq:scaling}, is a rough sub-Ricci flow if and only if $X$ is the unit sphere with the round metric and a multiple of volume measure.

\subsubsection*{Gaussian weights} 
Let $X=\R^n$.
Take a family of distances $d_t$ induced by the inner product $\langle \cdot, A_t\cdot\rangle$ where $A:I\to \R^{n\times n}$ is positive definite for all $t\in I$.
Consider $\m_t$ to be the weighted Lebesgue measure $e^{-f_t}\mathcal L^n$ with 
\[
    f_t(x)=\frac12\langle x,a_t x\rangle+\langle x,b_t\rangle+c_t\;, 
\]
where $a: I\to \R^{n\times n}$, $b: I\to \R^n$, $c:I\to\R$ are suitably regular functions.
Then $(X,d_t,\m_t)_{t\in I}$ is a weak/rough super-resp.~sub-Ricci flow if and only if 
\[\dot A_t \geq -2 a_t\;,\quad \text{resp.}\quad \dot A_t\leq -2a_t\;.\]
However, it will not be a $N$-super-Ricci flow for some $N\in [n,\infty)$ unless $a\equiv 0$ and $b\equiv 0$.

In particular, consider $A_t=1-2t$, $a_t\equiv 1$, $b_t\equiv 0$ and $c\equiv 0$ for $t\in(0,\frac12)$ corresponding to $d_t=\sqrt{1-2t}|\cdot|$ and $f_t(x)=\frac{|x|^2}{2}$.
Then the time-dependent space $(\R^n,d_t,e^{-f}\mathcal{L}^n)_{t\in(0,\frac12)}$, often referred to as a shrinking Gaussian, is a rough Ricci flow but not a non-collapsed Ricci flow.

\section{Minimal Super-Ricci Flows}\label{sec:minimal}
  
   In this part, we present an approach to synthetic characterisations of Ricci flows via the asymptotic behaviour of volumes. Let us illustrate the basic idea in the static case.
   \medskip
   
     In order to characterize Ricci flatness of closed Riemannian manifolds $(M,g)$, it suffices to characterize \emph{nonnegative Ricci curvature} and nonpositive scalar curvature or \emph{nonpositive integral of scalar curvature} since
     $$\Ric_x\ge0 \quad\textrm{and} \quad R(x)\le0\qquad\Longleftrightarrow\qquad \Ric_x=0$$
     for every $x\in M$, and 
      $$\Ric_x\ge0 \ (\forall x)\quad\textrm{and} \quad \int_M R\,d\vol\le0\qquad\Longleftrightarrow\qquad \Ric_x=0 \ (\forall x).$$
    For a synthetic characterization of scalar curvature and its integral, we propose to consider the  initial negative slope 
    $\k(x):=\liminf_{s\searrow0} \frac1s\big( \A_0(x)-\A_s(x)\big)$ of the 
      Gaussian volume functional
 $$\A_s(x):=(12\pi s)^{-n/2}\int_{M} e^{-\frac{d^2(x,y)}{12s}}m(dy)$$
 as well as the initial negative slope ${\mathfrak K}:=\liminf_{s\searrow0} \frac1s\big( {\mathfrak A}_0-{\mathfrak A}_s\big)$
of the Gaussian double integral
 ${\mathfrak A}_s:= \int_{M}\A_s(x)\, m(dx)$.

 \subsection{Ricci Bounds}
 
 \subsubsection*{The Gaussian Volume Functional and the Gaussian Double Integral}
 Let $(M,d,m)$ be an infinitesimally Hilbertian mm-space with finite volume and let $n\in\N$  be given.
   Assume that for every $x\in M$  the limit
 \begin{equation}\label{density}
 \rho(x):=\lim_{r\searrow0}\frac{m(B(r,x))}{\omega_nr^n}\end{equation}
 exists and is finite $m$-a.e.~where $\omega_n:=\frac{\pi^{n/2}}{\Gamma(n/2+1)}$ denotes the volume of the unit ball in $\R^n$.

The fundamental objects for the subsequent discussion are  the \emph{Gaussian volume functional} 
 $$\A_s(x):=(12\pi s)^{-n/2}\int_{M} e^{-\frac{d^2(x,y)}{12s}}m(dy),$$
its limit $\A_0(x):=\lim_{s\searrow0}\A_s(x)=\rho(x)$, and its initial negative slope  
$$\k(x):=\liminf_{s\searrow0} \frac1s\Big( \A_0(x)-\A_s(x)\Big),$$
as well as 
the \emph{Gaussian double integral} 
 $${\mathfrak A}_s:=(12\pi s)^{-n/2}\int_{M} \int_{M}e^{-\frac{d^2(x,y)}{12s}}m(dy)\, m(dx)$$
 its limit $\mathfrak A_0:=\lim_{s\searrow0}\mathfrak A_s=\int\rho\,dm$, and its initial negative slope 
 $${\mathfrak K}:=\liminf_{s\searrow0} \frac1s\Big( {\mathfrak A}_0-{\mathfrak A}_s\Big),$$
 called \emph{curvature functional}.
 Note that all these quantities depend on $n$. To emphasize this dependence,
 we occasionally write $\k^{(n)}$ and ${\mathfrak K}^{(n)}$.

The slope functional $\k$ provides a synthetic way to define the
scalar curvature.
For singular spaces, however, the slope functional $\mathfrak K$ turns out to provide more detailed information on the ``singular contributions to the curvature''. 
We will see that under minimal assumptions
${\mathfrak K}\ge \int_M \k\,dm$, in smooth cases ${\mathfrak K}= \int_M \k\,dm$,
but in numerous cases
${\mathfrak K}\not= \int_M \k\,dm$. 

 \begin{proposition}
  Let  $(M,g)$ be a smooth closed $n$-dimensional Riemannian manifold. Denote its
  scalar curvature by 
  $R=R_g$ and put $d:=d_g, m:=\vol_g$. Then
  $$\k(x)=R(x), \qquad
 {\mathfrak K}=\int_M R\,d\vol_g.$$
 \end{proposition}

 \subsubsection*{The Directional Gaussian Volume}
 A refined version of the Gaussian volume functional also provides a characterization of the Ricci curvature (rather than scalar curvature). To do so, requires to pass to the directional version. 
 
 Let $(M,d,m)$ be as before. In addition, now assume that
 \begin{itemize}
 \item it admits a unique tangent cone $T_xM$ at each point $x\in M$ and an associated unit tangent bundle $S_xM$ equipped with a Borel measure $\sigma_x$
 \item for each $x$ there exists a $\delta_x>0$ and an exponential map $\exp_x: 
 (0,\delta)\times S_xM\to M, (t,z)\mapsto \exp_x(tz)$.
 \end{itemize}
Given a measurable subset $Z\subset S_xM$, define the set 
 $$Z(x):=\Big\{ \exp_x(rz): z\in Z, r\in [0,\delta_x]\Big\},$$
and put
 $$\A_s(x, Z):=(12\pi s)^{-n/2}\int_{Z(x)} e^{-\frac{d^2(x,y)}{12s}}m(dy)$$
 as well as $$\k(x,Z):=\liminf_{s\to0} \frac1s\Big( \A_0(x, Z)-\A_s(x, Z)\Big).$$
 
%

 \begin{theorem}\label{k-ric} For a smooth (not necessarily complete) $n$-dimensional Riemannian manifold $(M,g)$ without boundary,
$$ \k(x,Z)=
\frac n{|S_x|}\int_{Z}\Ric(z)\,d\sigma(z).$$
 \end{theorem}

%
\begin{theorem} For a smooth complete $n$-dimensional Riemannian manifold $(M,g)$ without boundary and a number $K\in\R$ the following are equivalent:
\begin{enumerate}[(i)]
\item Lower Ricci bound: $\Ric_x\ge K g_x$ for all $x$
\item Bakry-\'Emery condition $\BE(K,n)$
\item Curvature-dimension condition $\CD(K,n)$
\item Measure contraction property $\MCP(K,n)$
\item Sharp Laplace comparison: for all $x$ and a.e. $y\in M$,
$$\frac12\Delta d^2(x,y) \le \begin{cases}1+(n-1)\frac{\sqrt{\frac K{n-1}}\,d(x,y)}{\tan\left(\sqrt{\frac K{n-1}}\,d(x,y)\right)} &\quad\text{if }K>0,\\
n &\quad\text{if }K=0,\\
1+(n-1)\frac{\sqrt{\frac{-K}{n-1}}\,d(x,y)}{\tanh\left(\sqrt{\frac{-K}{n-1}}\,d(x,y)\right)} &\quad\text{if }K<0.
\end{cases}
$$

\item Weak Laplace comparison: for all $x$ and a.e. $y\in M$,
$$\frac12 \Delta_y d^2(x,y)\le  n-\frac13 K d^2(x,y)$$
\item Perceived lower Ricci bound: for all $x$ and $Z$,
$$\k(x,Z)\ge {Kn}\,\frac{|Z|}{|S_x|} .$$
\end{enumerate}

\end{theorem}

  \subsection{Examples for Ricci bounds}
 
%
 
 \subsubsection*{Weighted Riemannian Manifolds}
Let $(M,g)$ be a smooth, closed $n$-dimensional Riemannian manifold  and put $d=d_g$,  $dm=\rho\,d\vol_g$ for some $\rho\in\mathcal C^2(M)$. The triple $(M,d,\rho)$ is called {weighted Riemannian manifold}. It obviously satisfies \eqref{density}.
We have
 $$\A_s(x)=(12\pi s)^{-n/2}\int e^{-\frac{d^2(x,y)}{12s}}{\rho(y)}d\vol_g(y), \qquad \A_0(x)=\rho(x),$$
 and
 $${\mathfrak A}_s=\int \A_s(x)\rho(x)\,d\vol_g(x), \qquad {\mathfrak A}_0=\int \rho^2d\vol.$$
 \begin{proposition}
$$ {\mathfrak k}(x)
=R(x)\rho(x)-3\Delta \rho (x),\qquad
 {\mathfrak K} =
\int \Big[\rho^2 R+  3|\nabla \rho|^2\Big]d\vol_g(x).
$$
\end{proposition}
 
 \begin{remark} Define the Bakry-\'Emery Ricci tensor by $\Ric^\rho=\Ric- \text{\rm Hess} \log\rho$ and the associated trace 
 by  $R^\rho=R-\Delta \log\rho$. 
 Then 
$$ {\mathfrak K} =\int \Big[\rho^2 R^\rho+  |\nabla \rho|^2\Big]d\vol_g(x).
$$

 \end{remark}
 

\begin{theorem} Assume  $(M,g, \rho)$ is a weighted $n$-dimensional Riemannian manifold as above with  $\Ric\ge 0$ or with $\Ric^\rho\ge 0$ and that it is \emph{minimal} in the sense that
${\mathfrak K}^{(n)}\le 0$. 
Then $\rho$ is constant and $M$ is \emph{Ricci flat}: $$\Ric_x=0 \quad\text{for all $x$.}$$
\end{theorem}

  \subsubsection*{RCD Spaces}
 Let $(M,d,m)$ be an infinitesimally Hilbertian mm-space which locally satisfies  the $\RCD(k,n)$-condition for some lower semicontinuous functions $k:M\to\R$ 
 (cf. \cite{BraunHabermannSturm,SturmDistrvaluedRicci,KettererVariablecurvbounds}) and such that 
 \eqref{density} holds. 
%
  \begin{theorem} Then
$$ {\mathfrak k}^{(n)}(x)
\ge nk(x)\rho(x).$$
 \end{theorem}

\begin{corollary} Assume in addition that $k$ is lower bounded and that $M$ has finite volume. Then
$$ {\mathfrak K}^{(n)}\ge n \int  k\rho\,dm.
$$
\end{corollary}
Note that 
if $(M,d,m)$  is non-collapsed (cf. \cite{DePhilippisGiglincRCD}), then $\rho\le1$ on $M$ and $\rho=1$ a.e. on $M$. 
\medskip

\subsubsection*{Manifolds with Boundary, Doubling, Tripling etc.} 
  \paragraph{A. Doubling} 
    \begin{theorem}  Let  $(M,g)$ be a compact Riemannian manifold with smooth boundary.  
Assume that $(\hat M, \hat g)$ is obtained by gluing two copies of $(M,g)$ along the boundary $\partial M$.
Then
 $$\widehat{\mathfrak K}
 =2\,\int_M R(x)d\vol(x)+4(n-1)\, \int_{\partial M} H(z)d\sigma(z)$$
 with $H(z)=$ mean curvature of $\partial M$ at $z\in\partial M$.
  \end{theorem}
 Note that $\int_{\hat M}\hat\k(x) dm(x)=2\,\int_M R(x)d\vol(x)$.
 
 \begin{example} Let $(M,g)$ be the closed ball of radius $\rho$ in $\R^n$ and let $(\hat M,\hat g)$ be its doubling which can be represented as warped product $[-\rho,\rho]\times_f {\mathbb S}^{n-1}$
with warping function
$f(r)=\rho-|r|$.
Approximate $(\hat M,\hat g)$ by Riemannian manifolds $(\hat M_\epsilon, \hat g_\epsilon)$ given as warped products $[-\rho,\rho]\times_{f_\epsilon} M$
with 
$$f_\epsilon(r)=\begin{cases}\rho-|r|, &\qquad |r|\in [\epsilon\frac\pi2,\rho]\\
\epsilon\cos\big(r/\epsilon\big)+\rho-\epsilon\frac\pi2, &\qquad |r|\in [0,\epsilon\frac\pi2].
\end{cases}$$
Let $\widehat{\mathfrak K}_\epsilon$ be the curvature functional for $(\hat M_\epsilon, \hat g_\epsilon)$. Then
$$\lim_{\epsilon\to0}\widehat{\mathfrak K}_\epsilon=
\lim_{\epsilon\to0}\int_{\hat M_\epsilon}R_\epsilon\,d\vol_\epsilon=
4n(n-1)\omega_n\rho^{n-2}=4(n-1)\, \int_{\partial M} H(z)d\sigma(z)=\widehat{\mathfrak K}.$$
 \end{example}

\begin{corollary} Suppose $(M,g)$ is a compact two-dimensional Riemannian manifold with boundary, and that $(\hat M_\epsilon, \hat g_\epsilon)$, $\epsilon>0$, is a family of Riemannian surfaces obtained from the doubling $(\hat M, \hat g)$ by smoothing the metric in the neighborhood of the rim without changing the topology. Then by the Gau{\ss}-Bonnet Theorem, the associated curvature functionals $\widehat{\mathfrak K}_\epsilon$ are independent of $\epsilon$ and
$$\widehat{\mathfrak K}=\widehat{\mathfrak K}_\epsilon=8\pi\, \chi(M)
$$
where  $\chi(M)$ denotes the Euler characteristic of $M$.
\end{corollary}  
 \paragraph{B. Manifolds with Boundary} 
 \begin{theorem} 
 Let  $M$ be the any compact subset  in $ \R^n$ with smooth boundary of non-vanishing surface measure. Then
 $${\mathfrak A}_s=|M|-\sqrt{\frac{3s}\pi}\cdot |\partial M|+{\mathcal O}(s).$$
Thus $\mathfrak K=+\infty$.  \end{theorem}
 Note, however, that
$\int \k(x)dx=0$.

  \paragraph{C. Gluing Multiple Copies}

   \begin{theorem}  Let  $M$ be the any compact subset  in $ \R^n$ with smooth boundary of non-vanishing surface measure.
Assume that $\hat M$ is obtained by gluing $k$ copies of $M$ along the boundary $\partial M$.
Then $$\widehat{\mathfrak A}_s=k |M|+(k-2)\sqrt{\frac{3s}\pi}\cdot |\partial M|+{\mathcal O}(s).$$
 Therefore,  $\widehat{\mathfrak K}=-\infty$ 
 whenever
 $k>2$.
  \end{theorem}
 
 \subsubsection*{Cones} 

 \paragraph{A. Two-dimensional Case}
  \begin{theorem}
   \label{2d-cone}
   Assume that $M$ is a closed surface, smooth outside some conical singularities $z_1, \ldots, z_k$ with conical angles $2\pi-\alpha_1,\ldots, 2\pi-\alpha_k$ for some $\alpha_i\in (0,\pi]$.
   \begin{itemize}
  \item[(i)]
  Then  
   $${\mathfrak K}=
   \int_{M\setminus \{z_1,\ldots,z_k\}} R\,d\vol+\sum_i C_{\alpha_i}$$
   with
   $$C_\alpha:=6\left(1-\frac\alpha{2\pi}\right)\, 
\frac{\alpha-\sin\alpha}{1-\cos\alpha}=2\alpha-\frac{\alpha^2}{\pi}+\frac{\alpha^3}{15}+{\mathcal O}(\alpha^4)\quad\textrm{as }\alpha\to0.$$

 \item[(ii)] If $(M^\epsilon,g^\epsilon)$ denotes a  smoothing of $(M,g)$ and $R^\epsilon$ its scalar curvature, then in general
   $${\mathfrak K}\not=\lim_{\epsilon\to0}{\mathfrak K}^\epsilon.$$
    In other words,
     the correction term $\sum_i C_{\alpha_i}$ does not coincide with 
   $$\lim_{\epsilon\to 0}\int_{M^\epsilon} R^\epsilon(x)d\vol(x)-\int_{M\setminus \{z_1,\ldots,z_k\}} R(x)d\vol(x)=\sum_i 2\alpha_i.$$

   \item[(iii)] In general, 
   ${\mathfrak K}\not=\int\k d\vol$
   and 
   $\k(z_i)=0$ for $i=1,\ldots,k$.
   \end{itemize}

   \end{theorem}

\paragraph{B. Higher-dimensional Case}
 \begin{theorem}
 Assume that $M$ is a closed Riemannian manifold of dimension $n>2$, smooth outside some singularities $z_1, \ldots, z_k$, and that  
 in the $\delta_i$-neighborhood of $z_i$ it is isometric to the cone over the $(n-1)$-sphere with radius $\rho_i$. 
 Then
 $${\mathfrak K}=\int_{M\setminus \{z_1,\ldots,z_k\}}  R\,d\vol=\int_M\k \,d\vol.$$
 Moreover, if  $(M^\epsilon, g^\epsilon)$ denotes a smoothing of the singularities of $(M,g)$ by means of smooth concave warping functions then
 $${\mathfrak K}=\lim_{\epsilon\to 0}{\mathfrak K}_\epsilon=
\lim_{\epsilon\to 0}\int_{M^\epsilon} R^\epsilon\, d\vol^\epsilon.$$
 \end{theorem}

 \begin{corollary}
 \begin{itemize}
 \item[(i)]
 The cone  over $M_1={\mathbb S}^2(1/\sqrt 3)\times {\mathbb S}^2(1/\sqrt 3)$ is ``Ricci flat'' in the sense that $\Ric_x\ge0$ for all $x$ and $\mathfrak K\le0$.
 
  \item[(ii)]
 The cone  over $M_1={\mathbb S}^{n-1}(1)/\Z_2$ for $n\ge 3$ is ``Ricci flat'' in the above sense.
 \end{itemize}
 \end{corollary}
Recall that Eguchi and Hanson \cite{EguchiHanson} constructed a smooth Ricci flat metric on $T^*S^2$ which provides an approximation for the latter space with $n=4$, cf.~\cite[pag.~11]{CheegerNotes2001} and \cite[Rmk.~5.4]{sturm2017remarks}.

  \subsection{Ricci Flows and Super-Ricci Flows}
  
  In the sequel, let $I$ denote a non-empty open subset of $\R$.
 
   \subsubsection*{The Gaussian Volume Functional and the Gaussian Double Integral}
Given a time-dependent family $(M,d_t,m_t)_{t\in I}$ of mm-spaces with $n$-dimensional volume densities 
 \begin{equation*}\label{densities}
 \rho_t(x):=\lim_{r\searrow0}\frac{m_t(B_t(r,x))}{\omega_nr^n},
 \end{equation*}
we define the \emph{time-dependent Gaussian volume functional} by
 $$\A_{s}(t,x):=(12\pi s)^{-n/2}\int e^{-\frac{d^2_{t+s}(x,y)}{12s}}m_t(dy),$$
  its limit  $\A_0(t,x):=\lim_{s\searrow0}\A_s(t,x)=\rho_t(x)$,  
 and \emph{its negative initial slope} by $$ \k(t,x):=\liminf_{s\searrow 0}\frac1s\Big[\A_{0}(t,x)-\A_{s}(t,x).
 \Big]$$
  Moreover, we define 
   \emph{the time-dependent Gaussian double integral}  by
  $${\mathfrak A}_t :=  (12\pi s)^{-n/2}\int\int e^{-\frac{d^2_{t+s}(x,y)}{12s}}\, m_t(dy)\,m_t(dx),
  $$
  its limit  ${\mathfrak A}_0:=\lim_{s\searrow0}{\mathfrak A}_t=\int\rho_t\,dm_t$,  
  and its initial negative slope by
  $${\mathfrak K}:=\liminf_{s\searrow0}\frac1s\left({\mathfrak A}_0-{\mathfrak A}_t
  \right).$$

\begin{theorem}  Assume that $(M,g_t)_{t\in I}$ is a smooth time-dependent Riemannian manifold and satisfies
 $\Ric_{t,x}+\frac12\partial_t g_{t,x} \ge0$  for given $t$ and $x$. Then
   $$\k(t,x) \ge0
  $$
  and 
  $$\k(t,x)=0 \quad \Longleftrightarrow \quad \Ric_{t,x}+\frac12\partial_t g_{t,x}=0.$$ 
 \end{theorem}

\begin{theorem}  Assume  that the smooth time-dependent Riemannian manifold $(M,g_t)_{t\in I}$ evolves as super-Ricci flow, that is, 
 $\Ric_{t}+\frac12\partial_t g_t \ge0\text{ on }I\times M$. Then for every $t\in I$,
   $${\mathfrak K}_t \ge0
  $$
  and 
  $${\mathfrak K}_t=0 \quad \Longleftrightarrow \quad \Ric_{t}+\frac12\partial_t g_t=0\text{ on }\{t\}\times M.$$ 
 \end{theorem}
 
%
%

\subsubsection*{Directional Gaussian Volume}
Similarly to the time-independent case, for given $t$ and $x$ we
define the set
 $$Z(t,x):=\Big\{ \exp_{t,x}(rz): z\in Z, r\in [0,\delta_x]\Big\}$$
 for any measurable subset $Z\subset S_{t,x}M$, and put
 $$\A_s(t,x, Z):=(12\pi s)^{-n/2}\int_{Z(t,x)} e^{-\frac{d_{t+s}^2(x,y)}{12s}}m_t(dy)$$
 and $$\k(t,x,Z):=\liminf_{s\to0} \frac1s\Big( \A_0(t,x, Z)-\A_s(t,x, Z)\Big).$$

 \begin{theorem}\label{ktZ}
 In the smooth Riemannian case,
$$\k(t,x,Z)=
\frac n{|S_{t,x}|}\int_{Z}\Big[\Ric(z,z)+\frac12\partial_tg_t(z,z)\Big]\,d\sigma_{t,x}(z).$$
 \end{theorem}

\paragraph{Novel Characterizations of SRFs}

\begin{theorem} For every smooth time-dependent Riemannian manifold $(M,g_t)_{t\in I}$, the following are equivalent:
\begin{itemize}
\item[(i)] For every $t$ and $x$, $$\Ric_{t,x}+\frac12\partial_t g_{t,x}\ge0.$$
\item[(ii)]
For every $t$ and every unit speed $d_t$-geodesic $\gamma$, 
$$\frac1{2}\Delta_y d_t^2(\gamma_0,\gamma_r) \le n+\frac2r\int_0^r \partial_{t}d_t(\gamma_q,\gamma_r)\,q\,dq.$$
\item[(iii)] For every $t$, $x$ and $Z$,
$$\k(t,x,Z)\ge0.$$
\end{itemize}
\end{theorem}

\subsection{Examples for (super-) Ricci flows}
Throughout the sequel,  let $I$ will denote a non-empty open subset of $\R$ and 
  $(M,d_t,m_t)_{t\in I}$ a time-dependent family of mm-spaces  with $n$-dimensional volume densities $\rho_t$ according to
 \eqref{densities}. Moreover, we assume that for all $s,t\in I$ with $s\le t$ the dual heat operator
 $\hat P_{t,s}: \mathcal P(M)\to \mathcal P(M)$
 exists.
 \begin{definition}\label{mSRF} (i) We say that  $(M,d_t,m_t)_{t\in I}$ is a \emph{super-Ricci flow (SRF in short)} if
 \begin{equation}\label{SRF}W_s\big(\hat P_{t,s}\mu,\hat P_{t,s}\nu\big)\le W_t\big(\mu,\nu\big)\end{equation}
 for all $\mu,\nu\in \mathcal P(M)$ and all $s,t\in I$ with $s\le t$.
 
 (ii) We say that  $(M,d_t,m_t)_{t\in I}$ is a \emph{minimal super-Ricci flow} if in addition for all $t\in I$,
$${\mathfrak K}_t\le0.$$
\end{definition}
Recall that in the smooth Riemannian case,  $(M,d_t,m_t)_{t\in I}$ is a minimal super-Ricci flow if and only if it is a Ricci flow.
Let us consider this concept of minimal SRF for a number of singular spaces. 

\begin{remark} In all the cases of time-dependent mm-spaces  $(M,d_t,m_t)_{t\in I}$ to be considered below, the definition \eqref{SRF} of SRF can equivalently be replaced by any other condition considered in \cite{Kopfer-Sturm2018,KopferSturm21}, including
\begin{itemize}
\item[(i)] Dynamic convexity of entropy 
$$
\partial_a \mathrm{Ent}(\mu^{a}|m_t)\big|_{a=1}-\partial_a \mathrm{Ent}(\mu^{a}|m_t)\big|_{a=0}
\ge- \frac 12\partial_tW_{t}^2(\mu^0,\mu^1)$$
\item[(ii)] Dynamic Bakry-\'Emery inequality 
$$\Gamma_{2,t}(f,f)\ge \frac12\partial_t \Gamma_{1,t}(f,f)$$

\item[(iii)]  
Gradient estimate
$$|\nabla_t P_{t,s}f|_t^2\le P_{t,s}|\nabla_s f|_s^2$$
\item[(iv)] 
Logarithmic Harnack inequality
\begin{align*}
 P_{t,s}(\log f)(x)  \leq \log (P_{t,s}f)(y)+\frac{d_t^2(x,y)}{4(t-s)}.
 \end{align*}
\end{itemize}

\end{remark}
\subsubsection*{Time-dependent Weighted Riemannian Manifolds}

 Let a family of smooth Riemannian manifolds $(M,g_t)_{t\in I}$ and a family  $(\rho_t)_{t\in I}$ of $\mathcal C^2$ functions on $M$ be given. Put $d_t=d_{g_t}$ and $m_t={\rho_t}\vol_t$.
 Then the Gausssian volume functional is 
 $$\A_{s}(t,x):=(12\pi s)^{-n/2}\int e^{-\frac{d^2_{t+s}(x,y)}{12s}}{\rho_t}\vol_t(dy)$$
 and, for instance,  the curvature functional, that is, the
initial negative slope of the time-dependent Gaussian double integral is
  $${\mathfrak K}_t :=\liminf_{s\searrow}\frac1s\left(\int\rho_t^2(x)\vol_t(dx)-  \int\int (12\pi s)^{-n/2}e^{-\frac{d^2_{t+s}(x,y)}{12s}}{\rho_t(y)}\vol_t(dy){\rho_t(x)}\vol_t(dx)\right).
  $$

  \begin{theorem}  (i) Assume 
 $\Ric_{t}+\frac12\partial_t g_t \, {\ge0}$. Then
 $$ \k(t,x) \ge-3\Delta {\rho_t}(x)$$
and
 $${\mathfrak K}_t
 \ge 3\int|\nabla\rho_t|^2d\vol_t.$$
(ii) Assume 
 $\Ric_{t}-\text{\rm Hess}_t\, \log\rho_t+\frac12\partial_t g_t \, \ge0$. Then
 $$\k(t,x) \ge-3\Delta {\rho_t}(x) +\rho_t\Delta\log\rho_t
 $$
 and
 $${\mathfrak K}_t 
 \ge \int|\nabla\rho_t|^2d\vol_t.$$
 \end{theorem}
  
%
%
%

\begin{corollary}  Assume that  $(M,g_t,\rho_t)$ satisfies either $\Ric_{g_t}+\frac12\partial_t g_t\ge0$ or  
 $\Ric_{t}-\text{Hess}_t\, \log\rho_t+\frac12\partial_t g_t \, \ge0$. Then for every $t\in I$,
   $${\mathfrak K}_t=0 \qquad\Longleftrightarrow\qquad \textrm{$\rho_t\equiv$ const in $x$  and $\Ric_{g_t}+\frac12\partial_t g_t=0$. }$$
 \end{corollary}

\begin{corollary}    $(M,d_{g_t},\rho_t\vol_{g_t})$  is a minimal SRF in the sense of Definition \ref{mSRF} if and only if $\rho_t\equiv$const ($\forall t$) and $(M,g_t)_{t\in I}$ is a Ricci flow. \end{corollary}

\subsubsection*{Surfaces of Convex Polytopes}
Let $(M,d,m)$ be the surface of an $(n+1)$-dimensional compact, convex polytope $\hat M\subset \R^{n+1}$. Put $d_t:=d, m_t:=m$ for all $t\in I$.
\begin{proposition} (i) The static family $(M,d_t,m_t)_{t\in I}$ is a SRF.

(ii) The SRF $(M,d_t,m_t)_{t\in I}$ is not minimal. 
\end{proposition}

\subsubsection*{Suspensions of Circles}
Let $N:=\R/{L \Z}\simeq[0,L)$ be a circle of length $L$ and $M:=[0,2\pi] \times_f N$  with $f(r):=\sin r$ the spherical suspension.
\begin{proposition} (i) If $L<2\pi$ then the linear shrinking space $(M,g_t)_{t\in I}$ with $I=(0,\frac12)$ and   $g_t=(1-2t)\,g$ is  a SRF but not minimal.

(ii) If $L>2\pi$ then there exist no rescaling $\alpha:(0,\delta)\to\R_+$ s.t. $(M,\alpha_tg)_{t\in (0,\delta)}$ is a SRF.

\end{proposition}

\subsubsection*{Suspensions in $n\ge3$}

\begin{proposition} Let $(M,g)$ be the spherical suspension over the round sphere $N:={\mathbb S}^{n-1}(\rho)$ of dimension $n-1$ and radius $\rho$.

(i) 
 If $\rho<1$ then the linear shrinking space $(M,g_t)_{t\in I}$ with $I=(0,\frac1{2(n-1)})$ and   $g_t=(1-2(n-1)t)\,g$ is  a SRF but not minimal.

(ii)  If $\rho>1$ then there exist no rescaling $\alpha:(0,\delta)\to\R_+$ s.t. $(M,\alpha_tg)_{t\in (0,\delta)}$ is a SRF.
\end{proposition}
\begin{proposition} Let $(M,g)$ be the spherical suspension over $(N,h)$ 
where 
 \begin{itemize}
 \item
 either $N={\mathbb S}^2(1/\sqrt 3)\times {\mathbb S}^2(1/\sqrt 3)$ and $n=5$
 \item or $N={\mathbb S}^{n-1}(1)/\Z_2$ for some $n\ge3$.
 \end{itemize}
 Note that $(M,g)$ is not a Riemannian manifold.
 
 Put $I=(0,\frac1{2(n-1)})$ and $g_t:=(1-2(n-1)t)g$. Then $(M,g_t)_{t\in I}$ 
 is a minimal SRF.
\end{proposition}

\subsubsection*{Subsets   $F\subset \R^n$ and their Doubling}
Let $F$ be a closed connected  set in $\R^n, n\ge2$, with smooth boundary and $\emptyset\not= F\not=\R^n$. Let $d$ be the induced euclidean distance and $m$ the Lebesgue measure on $F$.
\begin{proposition}
(i) If $F$ is convex then the static space $(F,d,m)$ is a SRF but not minimal.

(ii) If $F$ is nonconvex,  then the static space is not a SRF.
\end{proposition}

\begin{proposition} Let  $(\hat F,\hat d,\hat m)$ denote the doubling of  $(F,d,m)$ along $\partial F$.

(i) The static space  $(\hat F,\hat d,\hat m)$ is a SRF if and only if $F$ is convex.

(ii) The static space  $(\hat F,\hat d,\hat m)$ is a minimal SRF if and only if $F$ is a halfspace or a strip, that is, isometric either to $\R_+\times \R^{n-1}$ or to $[0,L]\times \R^{n-1}$ for some $L>0$.
\end{proposition}

\end{document}